\documentstyle[12pt]{amsart}
\newtheorem{theorem}{Theorem}[section]
\newtheorem{corollary}[theorem]{Corollary}
\newtheorem{remark}[theorem]{Remark}
\newtheorem{lemma}[theorem]{Lemma}
\newtheorem{proposition}[theorem]{Proposition}

\numberwithin{equation}{section}

\begin{document}

\title[ The mean curvature flow along the K\"ahler-Ricci flow]
{A mean curvature flow along a K\"ahler-Ricci flow}

\author{Xiaoli Han, Jiayu Li}

\address{Xiaoli Han, Department of Mathematical Sciences, Tsinghua University \\ Beijing 100084, P. R. of China.}
\email{xlhan@@math.tsinghua.edu.cn}

\address{Jiayu Li, Department of Mathematics, University of Science and Technology of China Hefei 230026 \\ AMSS CAS Beijing 100190, P. R. China}
\email{lijia@@amss.ac.cn}

\thanks {The first author was supported by NSF in China, No.10901088. The second author was supported by NSF in China, No. 11071236.}

\begin{abstract}
Let $(M,\overline{g})$ be a K\"ahler surface, and $\Sigma$ an immersed surface in $M$. The K\"ahler angle of $\Sigma$ in $M$ is introduced by
Chern-Wolfson \cite{CW}. Let $(M,\overline{g}(t))$ evolve along the K\"ahler-Ricci flow, and $\Sigma_t$ in $(M,\overline{g}(t))$ evolve along the
mean curvature flow. We show that the K\"ahler angle $\alpha(t)$ satisfies the evolution equation:
$$
(\frac{\partial}{\partial t}-\Delta)\cos\alpha=|\overline\nabla
J_{\Sigma_t}|^2\cos\alpha+R\sin^2\alpha\cos\alpha,
$$ where $R$ is the scalar curvature of $(M, \overline{g}(t))$.

The equation implies that, if the initial surface is symplectic (Lagrangian), then  along the flow, $\Sigma_t$ is
always symplectic (Lagrangian) at each time $t$, which we call a symplectic (Lagrangian) K\"ahler-Ricci mean curvature flow.

In this paper, we mainly study the symplectic K\"ahler-Ricci mean curvature flow.

\end{abstract}
\maketitle

\section{introduction }

Suppose that $(M, J, \omega, \overline g)$ is a K\"ahler surface.
Let $\Sigma$ be a compact oriented real surface which is smoothly
immersed in $M$, the K\"ahler angle $\alpha$ of $\Sigma$ in $M$ is
defined by Chern-Wolfson (\cite{CW})
\begin{equation}\label{al}
\omega|_\Sigma=\cos\alpha d\mu_\Sigma
\end{equation} where $d\mu_\Sigma$ is the area element of $\Sigma$
of the induced metric from $\overline{g}$. We say that $\Sigma$ is
a symplectic surface if $\cos\alpha>0$, $\Sigma$ is a holomorphic
curve if $\cos\alpha\equiv 1$ and $\Sigma$ is a Lagrangian surface
if $\cos\alpha\equiv 0$.

If $M$ is a K\"ahler-Einstein surface, a symplectic mean curvature
flow in $M$ was studied by Chen-Tian \cite{CT}, Chen-Li
\cite{CL1}, Wang \cite{Wa} and Han-Li \cite{HL}, etc. The main
point is that, along the mean curvature flow, the K\"ahler angle
satisfies a parabolic equation. However, without the Einstein
condition, one can not have the nice equation. In this paper we
find that, if $M$ evolves along the K\"ahler-Ricci flow, the
K\"ahler angle satisfies the same parabolic equation.

Now let $\overline{g}(t)$ evolve along the K\"ahler-Ricci flow on $M$
and $\Sigma$ evolve along the mean curvature flow in $(M,
\overline{g}(t))$, that is,
\begin{eqnarray}\label{main}
\left\{ \begin{array}{clcr} \frac{\partial }{\partial t}\overline{g}(t,\cdot) &= &-\overline{Ric}(g(t, \cdot))+\frac{r}{2}\overline{g}(t, \cdot)\\
\frac{d}{dt}F_t &= &\vec{H}\\ \overline {g}(0, \cdot) &= &\overline{g}_0\\
F(\cdot, 0) &= &F_0\end{array}\right.
\end{eqnarray} where $F_0: \Sigma_0\rightarrow (M, \overline{g}_0)$ is the
initial immersion and $\vec{H}$ is the mean curvature vector  of
$\Sigma_t$ in $(M, \overline{g}(t))$, and $r$ is a constant. We call it a K\"ahler-Ricci mean curvature flow,
denoted by $(M, \overline{g}(t),\Sigma_t)$.

The Ricci flow was introduced by Hamilton \cite{Ha} in order to
study the famous Poincar\'e conjecture, which was finally achieved
by Perelman (\cite{P1}, \cite{P2}, \cite{P3}). The K\"ahler-Ricci
flow was introduced by Cao \cite{C} to study Calabi conjecture,
which was studied by many authors (see \cite{CXT}, \cite{TZ},
\cite{ST1}, \cite{ST2}). The mean curvature flow was intensively
studied by Huisken \cite{H1}, \cite{H2}. Recall that, if
$C_1(M)<0$ ($=0,~>0$), choosing $r=-2$ ($0,~2$) and the initial
K\"ahler form with $c_1(M)$ as its K\"ahler class, Cao \cite{C}
proved that the K\"ahler-Ricci flow exists globally, and converges
to a K\"ahler-Einstein metric at infinity in the case that
$C_1(M)\leq 0$.

The main point of this paper is to derive the evolution equation
of the K\"ahler angle along the K\"ahler-Ricci mean curvature
flow. The purpose is to find symplectic minimal surface, and
especially holomorphic curves in K\"ahler surfaces.

Assume that $\Sigma_t$ evolves along the K\"ahler-Ricci mean curvature flow in $(M,
\overline{g}(t))$. We show that the evolution equation of $\cos\alpha$ is
$$
(\frac{\partial}{\partial t}-\Delta)\cos\alpha=|\overline\nabla
J_{\Sigma_t}|^2\cos\alpha+\overline{R}\sin^2\alpha\cos\alpha,
$$ where $\overline{R}$ is the scalar curvature of $(M, \overline{g}(t))$, and $|\overline\nabla
J_{\Sigma_t}|^2$ will be defined in (\ref{Jdefin}).

The same equation is obtained independently by Chen-Li \cite{CL1} and Wang \cite{Wa}
(also see \cite{CT}), for a symplectic mean curvature flow in a K\"ahler-Einstein surface. By the maximum principle
for parabolic equations, we see that, if the initial surface is symplectic (Lagrangian), then along the K\"ahler-Ricci mean curvature flow
$\Sigma_t$ is always symplectic (Lagrangian), which we call {\it a symplectic (Lagrangian) K\"ahler-Ricci mean curvature flow}. We will show
in a forthcoming paper \cite{HL2} that, Lagrangian is preserved by K\"ahler-Ricci mean curvature flow in any dimension.

In this paper, we will show that a symplectic K\"ahler-Ricci mean curvature flow does not develop Type I singularity
under the assumption that K\"ahler-Ricci flow does not develop any singularity.
If the K\"ahler surface is sufficiently close to a K\"ahler-Einstein surface and the initial surface is sufficiently close to a holomorphic curve,
then the symplectic K\"ahler-Ricci mean curvature flow exists globally and converges to a holomorphic curve in a K\"ahler-Einstein surface at infinity.

Suppose that $M=M_1\times M_2$ where $M_1, M_2$ are Riemann
surfaces with unit K\"ahler forms $\omega_1, \omega_2$. If $(M_1,
\overline{g}_1(0))$ and $(M_2, \overline{g}_2(0))$ have the same
average scalar curvature, and the initial surface is a graph with
$\langle e_1\times e_2,\omega_1\rangle>\frac{\sqrt{2}}{2}$ where
$e_1,e_2$ is an orthonormal frame of the initial surface, then the
K\"ahler-Ricci mean curvature flow exists globally. In addition,
if the scalar curvature of $M_1, M_2$ is positive, then the
K\"ahler-Ricci mean curvature flow converges to a totally geodesic
surface at infinity. The same result was also proved by
Chen-Li-Tian \cite{CLT} and Wang \cite{Wa} in the case that $M_1,
M_2$ have the same constant curvature.

Throughout this paper we will adopt the following ranges of
indices: $$A, B,\cdots=1, \cdots, 4,$$
 $$\alpha,\beta,\gamma, \cdots=3, 4,$$ $$i, j, k, \cdots=1,2.$$

\section{short time existence}
In this section, we show the short time existence for the
K\"ahler-Ricci mean curvature flow (\ref{main}).
\begin{theorem}\label{th1}
The evolution equation (\ref{main}) has a solution $(M, \overline{g}(t), \Sigma_t) $
for a short time with any smooth
compact initial surface $\Sigma_0$ at $t=0$, that is,
the solution of (\ref{main}) exists on a maximum time interval [0,
T).
\end{theorem}

{\it Proof.} It is well-known that there exists $T_1>0$ such that
the Ricci flow exists on $[0, T_1)$. Let $\Sigma$ evolves by the
mean curvature flow in $(M, \overline{g}(t))$ for $t<T_1$. Now we
use a trick of De Turck \cite{D} to prove the short time existence
of $g(t)$.

Recall the Gauss-Weingarten equation
\begin{eqnarray*}
\frac{\partial^2 F^\alpha}{\partial x^i\partial
x^j}-\Gamma^k_{ij}\frac{\partial F^\alpha}{\partial
x^k}+\overline{\Gamma}^\alpha_{\rho\sigma}\frac{\partial
F^\rho}{\partial x^i}\frac{\partial F^\sigma}{\partial
x^j}=h^\beta_{ij}v_\beta^\alpha,
\end{eqnarray*} $\Gamma^k_{ij}$ is the Christoffel symbol of $g(t)$
and $\overline{\Gamma}^\alpha_{\rho\sigma}$ is the Christoffel
symbol of $\overline{g}(t)$. Note that,
\begin{eqnarray*}
\Delta_{g(t)} F&=&g^{ij}\nabla_i\nabla_j F\\
&=&g^{ij}(\frac{\partial^2 F}{\partial x^i\partial
x^j}-\Gamma^k_{ij}\frac{\partial F}{\partial x^k}).
\end{eqnarray*} So the mean curvature flow equation can be written
as
\begin{eqnarray*}
\frac{\partial F^\alpha}{\partial t}=\Delta_{g(t)}
F^\alpha+g^{ij}\overline{\Gamma}^\alpha_{\rho\sigma}(t)\frac{\partial
F^\rho}{\partial x^i}\frac{\partial F^\sigma}{\partial x^j}
\end{eqnarray*} which is not the strictly parabolic equation. In
order to apply the standard theory of strictly parabolic equation
to get short time existence, we use a trick of De Turck by
modifying the flow through a diffeomorphism of the  parameter
space of $\Sigma$. We consider the following equation
\begin{eqnarray}\label{e4}
\frac{\partial \tilde{F}^\alpha}{\partial
t}=\Delta_{g(t)}\tilde{F}^\alpha+g^{ij}\overline{\Gamma}^\alpha_{\rho\sigma}(t)\frac{\partial
\tilde{F}^\rho}{\partial x^i}\frac{\partial
\tilde{F}^\sigma}{\partial x^j}+v^k\nabla_k\tilde{F}^\alpha,
\end{eqnarray} where the vector field $v^k$ will be chosen to
make the equation strictly parabolic. In fact, if $\tilde{F}(y,
t)$ is a solution of (\ref{e4}), then
\begin{eqnarray*}
F(x, t)=\tilde{F}(y(x, t), t)
\end{eqnarray*} satisfies  the equation
\begin{eqnarray*}
\frac{\partial F^\alpha}{\partial t} &=&\frac{\partial\tilde{F}^\alpha}{\partial
t}+\nabla_k\tilde{F}^\alpha\cdot\frac{dy^k}{dt}\\
&=&\Delta_{g(t)}\tilde{F}^\alpha+g^{ij}\overline{\Gamma}^\alpha_{\rho\sigma}(t)\frac{\partial
\tilde{F}^\rho}{\partial x^i}\frac{\partial\tilde{F}^\sigma}{\partial x^j}+(v^k+\frac{dy^k}{dt})\nabla_k\tilde{F}^\alpha\\
&=&\Delta_{g(t)}
F^\alpha+g^{ij}\overline{\Gamma}^\alpha_{\rho\sigma}(t)\frac{\partial
F^\rho}{\partial x^i}\frac{\partial F^\sigma}{\partial x^j}
\end{eqnarray*} by choosing
\begin{equation}
\left\{\begin{array}{clcr}\frac{dy^k}{dt} &= &-v^k(x, t),\\
y^k|_{t=0} &= &x^k.\end{array}\right.
\end{equation} Now we pick
\begin{eqnarray*}
v^k=g^{ij}(\Gamma^k_{ij}(t)-\Gamma^k_{ij}(0)).
\end{eqnarray*} The equation (\ref{e4}) becomes
\begin{eqnarray*}
\frac{\partial \tilde{F}^\alpha}{\partial
t}=g^{ij}\{\frac{\partial^2\tilde{F}^\alpha}{\partial x^i\partial
x^j}-\Gamma^k_{ij}(0)\frac{\partial\tilde{F}^\alpha}{\partial
x^k}\}+g^{ij}\overline{\Gamma}^\alpha_{\rho\sigma}(t)\frac{\partial
\tilde{F}^\rho}{\partial x^i}\frac{\partial\tilde{
F}^\sigma}{\partial x^j}
\end{eqnarray*} which is a strictly parabolic equation. Since $\overline{\Gamma}^\alpha_{\rho\sigma}(t)$ is
uniformly bounded in $[0, T)$, then by the standard theory of
parabolic equations we get the short time existence of the mean
curvature flow (c.f. \cite{LSU}). Thus there exists $T>0$ such
that the K\"ahler-Ricci mean curvature flow exists on [0, T).
\hfill Q. E. D.

If the K\"ahler-Ricci mean curvature flow blows up at $ T$, there
are two possibilities. One is the K\"ahler-Ricci flow blows up at
$T$, another one is the mean curvature flow blows up at $T$. Let's
state some fundamental results regarding the singularity of the
Ricci flow and the mean curvature flow.
\begin{theorem}\label{th3}\cite{S}
If the Ricci curvature is uniformly bounded under the Ricci flow
$\frac{\partial}{\partial t}g_{ij}=-2R_{ij}$ for all times $t\in
[0, T)$, then the solution can be extended beyond $T$.
\end{theorem}

\begin{theorem}\label{th5}\cite{H1}
If the second fundamental form is uniformly bounded under the mean
curvature flow $\frac{dF}{dt}=\vec{H}$ for all times $t\in [0,
T)$, then the solution can be extended beyond $T$.
\end{theorem}

\section{evolution equations}
In this section the evolution equations of the metric and the
second fundamental form of $\Sigma_t$ will be derived along the
K\"ahler-Ricci mean curvature flow (\ref{main}). In terms of
coordinates $\{x^i\}$ on $\Sigma_t$ and coordinates $\{y^A\}$ on
$(M,\overline{g}(t))$, the metric of $\Sigma$ can be expressed as
follows:
$$g_{ij}(x, t)=\overline{g}_{AB}(F(x, t), t)\frac{\partial F_t^A}{\partial x^i}\frac{\partial F_t^B}{\partial x^j}.$$
\allowdisplaybreaks

\begin{lemma}
The metric of $\Sigma_t$ satisfies the evolution equation
\begin{equation}\label{g}
\frac{\partial}{\partial t}g_{ij}=-2H^\alpha
h^\alpha_{ij}-\overline{R}_{ij}+\frac{r}{2}g_{ij}.
\end{equation}
\end{lemma}

{\it Proof.} It is clear that
\begin{eqnarray}\label{e9}
\frac{\partial}{\partial t}g_{ij} &=&\overline{g}_{AB,C}\frac{\partial F^C_t}{
\partial t}\frac{\partial F_t^A}{\partial
x^i}\frac{\partial F_t^B}{\partial
x^j}+\frac{\partial}{\partial
t}\overline{g}_{AB}\frac{\partial F_t^A}{\partial
x^i}\frac{\partial F_t^B}{\partial
x^j}+2\overline{g}_{AB}(t)\frac{\partial }{\partial
x^i}(\frac{\partial F_t^A}{\partial t})\frac{\partial
F_t^B}{\partial x^j}\nonumber\\ &=&-\overline{Ric}(\frac{\partial F}{\partial
x^i}, \frac{\partial F}{\partial
x^j})+\frac{r}{2}\overline{g}_{AB}\frac{\partial F_t^A}{\partial
x^i}\frac{\partial F_t^B}{\partial x^j}-2H^\alpha
h^\alpha_{ij}\nonumber\\ &=&-\overline{R}_{ij}+\frac{r}{2}g_{ij}-2H^\alpha
h^\alpha_{ij}\nonumber.
\end{eqnarray} \hfill Q. E. D.

From the evolution equation of the metric we can get the evolution
equation of the area element.
\begin{corollary}
Set $\tilde{R}=\frac{1}{2}g^{ij}\overline{R}_{ij}$, the area element of $\Sigma_t$ satisfies the following equation.
\begin{eqnarray}\label{e100}
\frac{d}{dt}d\mu_t=(-|H|^2-\tilde{R}+r)d\mu_t,
\end{eqnarray} and consequently,
\begin{equation}\label{e17}
\frac{d}{dt}(e^{-rt}\int_{\Sigma_t}d\mu_t)=-e^{-rt}\int_{\Sigma_t}(|H|^2+\tilde{R})d\mu_t.
\end{equation}
\end{corollary}

\begin{lemma}\label{A}
Under the flow (\ref{main}), the second fundamental form
$h^\alpha_{ij}$ satisfies the following equation.
\begin{eqnarray}\label{e5}
\frac{d}{dt}h^\alpha_{ij} &=&\Delta
h^\alpha_{ij}+\overline\nabla_k \overline{R}_{\alpha
ijk}+\overline\nabla_j \overline{R}_{\alpha kik}\nonumber\\
&&-2\overline{R}_{lijk}h^\alpha_{lk}+2\overline{R}_{\alpha\beta
jk}h^\beta_{ik}+2\overline{R}_{\alpha\beta
ik}h^\beta_{jk}\nonumber\\&&-\overline{R}_{lkik}h^\alpha_{lj}-\overline{R}_{lkjk}h^\alpha_{li}+\overline{R}_{\alpha
k\beta k}h^\beta_{ij}\nonumber\\ &&-h^\alpha_{im}(H^\gamma
h^\gamma_{mj}-h^\gamma_{mk}h^\gamma_{jk})\nonumber\\&&-h^\alpha_{mk}(h^\gamma_{mj}
h^\gamma_{ik}-h^\gamma_{mk}h^\gamma_{ij})\\&&-h^\beta_{ik}(h^\beta_{lj}
h^\alpha_{lk}-h^\beta_{lk}h^\alpha_{lj})\nonumber\\&&-h^\alpha_{jk}h^\beta_{ik}H^\beta+h^\beta_{ij}b^\alpha_\beta\nonumber\\
&&-\overline{R}_{\alpha\beta}h^\beta_{ij}+\frac{r}{2}h^\alpha_{ij}
-\frac{1}{2}(\overline{\nabla}_i \overline{R}_{j\alpha}+\overline{\nabla}_j \overline{R}_{i\alpha}-\overline{\nabla}_\alpha \overline{R}_{ij}).
\end{eqnarray} where
$\overline{\nabla}$ is the covariant derivative of $(M, \overline{g}(t))$ and $b^\alpha_\beta=\langle\frac{\partial}{\partial t} e_\alpha, e_\beta\rangle$. In
particular, $|A|^2$ satisfies the following equation along the
flow (\ref{main}).
\begin{eqnarray}\label{e6}
\frac{d}{dt}|A|^2&=&\Delta |A|^2-2|\nabla
A|^2+2[\overline\nabla_k \overline{R}_{\alpha ijk}+\overline\nabla_j
\overline{R}_{\alpha
kik}]h^\alpha_{ij}\nonumber\\&&-4\overline{R}_{lijk}h^\alpha_{lk}h^\alpha_{ij}+8\overline{R}_{\alpha\beta
jk}h^\beta_{ik}h^\alpha_{ij}-4\overline{R}_{lk
ik}h^\alpha_{lj}h^\alpha_{ij}+2\overline{R}_{\alpha k\beta
k}h^\beta_{ij}h^\alpha_{ij}\nonumber\\&&+2\sum_{\alpha,\gamma, i,
m}(\sum_{k}h^\alpha_{ik}h^\gamma_{mk}-h^\alpha_{mk}h^\gamma_{ik})^2+2\sum_{i,
j, m, k}(\sum_{\alpha}h^\alpha_{ij}h^\alpha_{mk})^2\nonumber\\&&
+2\overline{R}_{ik}h^\alpha_{ij}h^\alpha_{kj}-2\overline{R}_{\alpha\beta}h^\beta_{ij}h^\alpha_{ij}-r|A|^2\nonumber \\&&
-\frac{1}{2}h^\alpha_{ij}(\overline{\nabla}_i \overline{R}_{j\alpha}+\overline{\nabla}_j \overline{R}_{i\alpha}-\overline{\nabla}_\alpha \overline{R}_{ij}).
\end{eqnarray}
More generally, we have
\begin{eqnarray}\label{e71}
\frac{d}{dt}|\nabla^m A|^2&=&\triangle |\nabla^m
A|^2-2|\nabla^{m+1} A|^2+\sum_{i+j+k=m}\nabla^i A*\nabla^j
A*\nabla^k A*\nabla^m A,
\end{eqnarray} where we denote by $\nabla^i A*\nabla^j
A*\nabla^k A$ the linear combination of $\nabla^i A, \nabla^j A$
and $\nabla^k A$.
\end{lemma}

{\it Proof.}
By $(7.4)$ in \cite{Wa}, the Laplacian of $h^\alpha_{ij}$ satisfies
\begin{eqnarray}\label{e7}
\Delta h^\alpha_{ij} &=&H^\alpha_{,ij }-(\overline\nabla_k
\overline{R})_{\alpha ijk}-(\overline\nabla_j \overline{R})_{\alpha
ijk}\nonumber\\&&+2\overline{R}_{lijk}h^\alpha_{lk}-2\overline{R}_{\alpha\beta
jk}h^\beta_{ik}-2\overline{R}_{\alpha\beta ik}h^\beta_{jk}-\overline{R}_{\alpha
ij\beta}H^\beta\nonumber\\&&+\overline{R}_{lkik}h^\alpha_{lj}+\overline{R}_{lkjk}h^\alpha_{li}-\overline{R}_{\alpha
k\beta k}h^\beta_{ij}\nonumber\\&&+h^\alpha_{im}(H^\gamma
h^\gamma_{mj}-h^\gamma_{mk}h^\gamma_{jk})\nonumber\\&&+h^\alpha_{mk}(h^\gamma_{mj}h^\gamma_{ik}-h^\gamma_{mk}h^\gamma_{ij})\nonumber\\
&&+h^\beta_{ik}(h^\beta_{lj}h^\alpha_{lk}-h^\beta_{lk}h^\alpha_{lj}).
\end{eqnarray} Now we compute
$\frac{\partial}{\partial t}h^\alpha_{ij}$. Since
$h^\alpha_{ij}=\langle\overline\nabla_i e_j,
e_\alpha\rangle_{\bar{g}}=\overline{g}_{\alpha A}(\overline\nabla_i
e_j)^A e_\alpha$, we get
\begin{eqnarray}\label{e8}
\frac{d}{dt}h^\alpha_{ij} &=&\frac{d}{dt}\overline{g}_{\alpha
A}(\overline\nabla_i e_j)^A,
e_\alpha+\langle\frac{\partial}{\partial t}(\overline\nabla) (e_i e_j),
e_\alpha\rangle+\langle\overline{\nabla}_H\overline{\nabla}_i e_j, e_\alpha\rangle+\langle\overline{\nabla}_i e_j, \frac{\partial}{\partial t}e_\alpha\rangle.\nonumber\\
\end{eqnarray} Using $(7.5)$ in \cite{Wa}, we have
\begin{eqnarray*}
\langle\overline\nabla_H\overline{\nabla}_i e_j, e_\alpha\rangle=H^\alpha_{,ij}-H^\beta
h^\beta_{ik}h^\alpha_{jk}-H^\beta \overline{R}_{\beta
ji\alpha}.
\end{eqnarray*}
Using $(2.24)$ in \cite{CLN}, we get that
\begin{eqnarray*}
\langle(\frac{\partial}{\partial t}\overline\nabla)(e_i, e_j)
, e_\alpha\rangle=-\frac{1}{2}(\overline{\nabla}_i \overline{R}_{j\alpha}+\overline{\nabla}_j \overline{R}_{i\alpha}-\overline{\nabla}_\alpha \overline{R}_{ij}).
\end{eqnarray*} Set $\langle\frac{\partial}{\partial t}e_\alpha, e_\beta\rangle=b^\beta_\alpha.$ Putting these equations into (\ref{e8}) we obtain that
\begin{eqnarray}\label{e97}
\frac{d}{dt}h^\alpha_{ij} &=& -(\overline{R}_{\alpha
\beta}-\frac{r}{2}\overline{g}_{\alpha\beta})
h^\beta_{ij}-\frac{1}{2}(\overline{\nabla}_i \overline{R}_{j\alpha}+\overline{\nabla}_j \overline{R}_{i\alpha}-\overline{\nabla}_\alpha \overline{R}_{ij})
\nonumber\\&&+H^\alpha_{,ij}-H^\beta
h^\beta_{ik}h^\alpha_{jk}-H^\beta R_{\beta
ji\alpha}+h^\beta_{ij}b^\alpha_\beta.
\end{eqnarray}

Combining equations (\ref{e7}) and (\ref{e97}), we
get the parabolic equation (\ref{e5}) for $h^\alpha_{ij}$.

Since $|A|^2=g^{ik}g^{jl}h^\alpha_{ij}h^\alpha_{kl}$, by (\ref{g})
we have,
\begin{eqnarray*}
\frac{d}{dt}|A|^2
&=&2(\frac{d}{dt}g^{ik})h^\alpha_{ij}h^\alpha_{kj}+2(\frac{d}{dt}h^\alpha_{ij})h^\alpha_{ij}\\
&=&2(2H^\beta
h^\beta_{ik}+\overline{R}_{ik}-rg_{ik})h^\alpha_{ij}h^\alpha_{kj}\\&&
+2h^\alpha_{ij}[\Delta h^\alpha_{ij}+\overline\nabla_k \overline{R}_{\alpha
ijk}+\overline\nabla_j \overline{R}_{\alpha kik}\\
&&-2\overline{R}_{lijk}h^\alpha_{lk}+2\overline{R}_{\alpha\beta
jk}h^\beta_{ik}+2\overline{R}_{\alpha\beta
ik}h^\beta_{jk}\\&&-\overline{R}_{lkik}h^\alpha_{lj}-\overline{R}_{lkjk}h^\alpha_{li}+\overline{R}_{\alpha
k\beta k}h^\beta_{ij}\\ &&-h^\alpha_{im}(H^\gamma
h^\gamma_{mj}-h^\gamma_{mk}h^\gamma_{jk})\\&&-h^\alpha_{mk}(h^\gamma_{mj}
h^\gamma_{ik}-h^\gamma_{mk}h^\gamma_{ij})\\&&-h^\beta_{ik}(h^\beta_{lj}
h^\alpha_{lk}-h^\beta_{lk}h^\alpha_{lj})\\&&-h^\alpha_{jk}h^\beta_{ik}H^\beta+h^\beta_{ij}b^\alpha_\beta\\
&&-\overline{R}_{\alpha\beta}h^\beta_{ij}+\frac{r}{2}h^\alpha_{ij}
-\frac{1}{2}(\overline{\nabla}_i \overline{R}_{j\alpha}+\overline{\nabla}_j \overline{R}_{i\alpha}-\overline{\nabla}_\alpha \overline{R}_{ij})].
\end{eqnarray*} Using
\begin{eqnarray*}
\Delta (h^\alpha_{ij})^2=2|\nabla A|^2+2h^\alpha_{ij}\Delta
h^\alpha_{ij},
\end{eqnarray*} and the antisymmetric of $b^\alpha_\beta, \langle e_\beta, \overline{\nabla}_H e_\alpha\rangle_{\bar{g}}$,  we get that
\begin{eqnarray*}
\frac{d}{dt}|A|^2 &=&\Delta |A|^2-2|\nabla A|^2\\
&&+2\overline{R}_{ik}h^\alpha_{ij}h^\alpha_{kj}-2\overline{R}_{\alpha\beta}h^\alpha_{ij}
h^\beta_{ij}-r|A|^2 \\&&+2h^\alpha_{ij}[\overline\nabla_k
\overline{R}_{\alpha
ijk}+\overline\nabla_j \overline{R}_{\alpha kik}\\&&
-\frac{1}{2}(\overline{\nabla}_i \overline{R}_{j\alpha}+\overline{\nabla}_j \overline{R}_{i\alpha}-\overline{\nabla}_\alpha \overline{R}_{ij})\\
&&-2\overline{R}_{lijk}h^\alpha_{lk}+2\overline{R}_{\alpha\beta
jk}h^\beta_{ik}+2\overline{R}_{\alpha\beta
ik}h^\beta_{jk}\\&&-\overline{R}_{lkik}h^\alpha_{lj}-\overline{R}_{lkjk}h^\alpha_{li}+\overline{R}_{\alpha
k\beta
k}h^\beta_{ij}\\&&+h^\alpha_{im}h^\gamma_{mk}h^\gamma_{jk}-h^\alpha_{mk}(h^\gamma_{mj}h^\gamma_{ik}-h^\gamma_{mk}h^\gamma_{ij})
-h^\beta_{ik}(h^\beta_{lj}h^\alpha_{lk}-h^\beta_{lk}h^\alpha_{lj})].
\end{eqnarray*} The fourth order terms can be calculated as in
\cite{Wa}. The equation (\ref{e71}) can be proved similarly.

\hfill Q. E. D.

The following theorem follows easily from Theorem \ref{th3},
Theorem \ref{th5} and (\ref{e71}).
\begin{theorem}
If the Ricci curvature of $\overline{g}(t)$ is uniformly bounded,
and the second fundamental form of $\Sigma_t$ is uniformly bounded
under the K\"ahler-Ricci mean curvature flow for all time $t\in
[0, T)$, then the solution can be extended beyond $T$.
\end{theorem}

\section{The evolution of the K\"ahler angle along the flow}
This is the main section of this paper, in which we will derive the evolution equation of the
K\"ahler angle along the
K\"ahler-Ricci mean curvature flow.

Let $(M, \overline{g})$ be a K\"ahler surface. Let $\overline{g}(t)$
evolve along the K\"ahler-Ricci flow on $M$, and $\Sigma_t$ evolve along the mean curvature flow in $(M,
\overline{g}(t))$. Choose an orthonormal basis $\{e_1, e_2, e_3, e_4\}$
on $(M,\overline{g}(t))$ along $\Sigma_t$ such that $\{e_1, e_2\}$ is the basis of
$\Sigma_t$. Let $J_{\Sigma_t}$ be an almost complex structure in a tubular neighborhood
of $\Sigma_t$ on $(M,\overline{g}(t))$ with
\begin{equation}\label{Jdefin}
\left\{\begin{array}{clcr} J_{\Sigma_t}e_1&=&e_2\\
J_{\Sigma_t}e_2&=&-e_1\\ J_{\Sigma_t}e_3&=&e_4\\
J_{\Sigma_t}e_4&=&-e_3.
\end{array}\right.
\end{equation}

It is proved in \cite{CL1} that
\begin{eqnarray}\label{e11}
|\overline{\nabla}J_{\Sigma_t}|^2&=& |h_{1k}^4+h_{2k}^3|^2
+|h_{2k}^4-h_{1k}^3|^2\nonumber\\ &\geq&\frac{1}{2}|H|^2.
\end{eqnarray}

\begin{theorem}\label{th2}
Let $(M, \overline{g}(t),\Sigma_t)$ be a K\"ahler-Ricci mean curvature flow.
Then the evolution equation of $\cos\alpha$ is
\begin{equation}\label{cos}
(\frac{\partial}{\partial t}-\Delta)\cos\alpha=|\overline\nabla
J_{\Sigma_t}|^2\cos\alpha+\overline{R}\sin^2\alpha\cos\alpha,
\end{equation} where $\overline{R}$ is the scalar curvature of $(M, \overline{g}(t))$.
As a corollary, if the initial surface $\Sigma_0$ is symplectic, then along the flow, at each time
$t$, $\Sigma_t$ is symplectic, and if the initial surface $\Sigma_0$ is Lagrangian, then along the flow, at each time
$t$, $\Sigma_t$ is Lagrangian.
\end{theorem}

{\it Proof.} Choose an orthonormal basis $\{e_1, e_2, e_3, e_4\}$
on $(M, \overline{g}(t))$ along $\Sigma_t$ such that $\{e_1, e_2\}$ is the basis of
 $\Sigma_t$ and $\omega_t$ takes the form

\begin{equation}\label{e31}
\omega_t=\cos\alpha u_1\wedge\ u_2 + \cos \alpha u_3\wedge u_4
+\sin\alpha u_1\wedge u_3 - \sin\alpha u_2\wedge u_4,
\end{equation} where $\{u_1, u_2, u_3, u_4\}$ is the dual basis of $\{e_1, e_2, e_3, e_4\}$.

Using the evolution equation of the metric (\ref{g}) and
(\ref{e31}), we get that
\begin{eqnarray*}
\frac{\partial}{\partial t}\cos\alpha &=&\frac{\partial}{\partial
t}\frac{\omega(e_1,
e_2)}{\sqrt{\det(g_{ij})}}=\frac{\partial}{\partial
t}\frac{\langle Je_1, e_2\rangle_{\overline{g}}}{\sqrt{\det(g_{ij})}}=
\frac{\partial}{\partial
t}\frac{\overline{g}_{2A}(Je_1)^A e_2}{\sqrt{\det(g_{ij})}}\\
&=&\frac{-\overline{Ric}(Je_1, e_2)+r/2\omega(e_1, e_2)+\omega(\nabla_1 H,
e_2)+\omega(e_1, \nabla_2 H)}{\sqrt{\det(g_{ij})}}\\
&&-\frac{1}{2}\cos\alpha g^{ij}\frac{\partial}{\partial
t}g_{ij}\\&=& -\overline{Ric}(Je_1,
e_2)+\frac{r}{2}\cos\alpha+\sin\alpha(H^4,1+H^3,2)-|H|^2\cos\alpha\\
&&+\frac{1}{2}\cos\alpha(\overline{R}_{11}+\overline{R}_{22})-\frac{r}{2}\cos\alpha+|H^2|\cos\alpha\\
&=&-\overline{Ric}(Je_1,
e_2)+\sin\alpha(H^4,1+H^3,2)+\frac{1}{2}\cos\alpha(\overline{R}_{11}+\overline{R}_{22}).
\end{eqnarray*}

Recall the equation in Proposition $3.1$ and Lemma $3.2$ in
\cite{HL1} for $\cos\alpha$,
\begin{eqnarray*}
\Delta\cos\alpha
&=&-|\overline{\nabla}J_{\Sigma_t}|^2\cos\alpha+\sin\alpha(H^4,1+H^3,2)-\sin^2\alpha
\overline{Ric}(Je_1, e_2).
\end{eqnarray*} Thus we have
\begin{eqnarray*}
(\frac{\partial}{\partial t}-\Delta)\cos\alpha
&=&|\overline{\nabla}J_{\Sigma_t}|^2\cos\alpha-\cos^2\alpha \overline{Ric}(Je_1,
e_2)+\frac{1}{2}\cos\alpha(\overline{R}_{11}+\overline{R}_{22}).
\end{eqnarray*} Using the equation in Lemma $3.2$ in \cite{HL1}, $\overline{Ric}(Je_1,
e_2)=\frac{1}{\cos\alpha}(\overline{R}_{1212}+\overline{R}_{1234})$ we have
\begin{eqnarray*}
-\cos^2\alpha \overline{Ric}(Je_1,
e_2)+\frac{1}{2}\cos\alpha(\overline{R}_{11}+\overline{R}_{22})&=&
-\cos\alpha(\overline{R}_{1212}+\overline{R}_{1234})\\&&+\frac{1}{2}\cos\alpha(2\overline{R}_{1212}+\overline{R}_{3i3i}+\overline{R}_{4i4i})\\
&=&\frac{1}{2}\cos\alpha(\overline{R}_{3i3i}+\overline{R}_{4i4i}-2\overline{R}_{1234})\\&=&\overline{R}\cos\alpha\sin^2\alpha,
\end{eqnarray*} where the last equality was derived in Lemma $3.2$
of \cite{MW}. Therefore we proved the theorem. \hfill Q. E. D.

Then by the parabolic minimum principle, we obtain that
\begin{theorem}\label{th6}
Suppose the smooth solution of (\ref{main})
exists on $[0, T)$. Let $(M, \overline{g}(0))$ be a K\"ahler surface with nonnegative
scalar curvature.  If $\cos\alpha(x, 0)\geq c_0>0$, then
$$\cos\alpha(x, t)\geq c_0,$$ for all $t\in [0, T)$.
\end{theorem}

{\it Proof.} Recall the evolution equation of the scalar curvature
of $M$ under the K\"ahler-Ricci flow,
\begin{eqnarray}\label{e101}
\frac{\partial}{\partial t}\overline{R}=\frac{1}{2}\Delta \overline{R}+|\overline{Ric}|^2-\frac{r\overline{R}}{2}.
\end{eqnarray} Thus by the parabolic minimum principle, if the
scalar curvature of the initial surface is nonnegative, then the
scalar curvature of $\overline{R}(t)$ is nonnegative for all $t\in [0, T)$.
Using the parabolic minimum principle again and (\ref{cos}), we
proved the theorem. \hfill Q. E. D.

\begin{corollary}
 We can rewrite the equation (\ref{cos}) as
\begin{eqnarray}\label{sin}
(\frac{\partial}{\partial
t}-\Delta)\sin^2(\alpha/2)=-|\overline{\nabla}J_{\Sigma_t}|^2\cos\alpha-4\overline{R}\sin^2(\alpha/2)\cos^2(\alpha/2)\cos\alpha,
\end{eqnarray} which yields that, along the flow (\ref{main}),
the K\"ahler angle decrease to zero exponentially if
$\overline{R}\geq\delta>0$ and $\cos\alpha(\cdot, 0)\geq c_0$.
\end{corollary}

\begin{remark}
The same evolution equation for $\cos\alpha$ along the mean curvature flow in the case that $M$
is K\"ahler-Einstein surface was proved by Chen-Li \cite{CL1} and Wang \cite{Wa}.
\end{remark}

\begin{lemma}\label{l2}
$$|\nabla\alpha|^2\leq |\overline{\nabla} J_{\Sigma_t}|^2.$$
\end{lemma}

{\it Proof.} Using the frame in Theorem \ref{th2} and (\ref{e31}), it is easy to see that
\begin{eqnarray}\label{e28}
{\overline{\nabla}}_1\cos\alpha &=& \omega({\overline{\nabla}}_{e_1}e_1, e_2)+\omega(e_1,{\overline{\nabla}}_{e_1}e_2)\nonumber\\
&=&(h^4_{11}+h^3_{12})\sin\alpha.
\end{eqnarray}Similarly, we can get that
\begin{eqnarray}\label{e29}
{\overline{\nabla}}_2\cos\alpha &=& (h^3_{22}+h^4_{12})\sin\alpha.
\end{eqnarray} Therefore,
$$|\nabla\alpha|^2\leq |\overline{\nabla} J_{\Sigma_t}|^2.$$ \hfill Q. E. D.

\section{Monotonicity formula}

In this section we assume that the K\"ahler-Ricci flow exists on
$[0, T]$, i.e, the Ricci curvature of $(M, \overline{g}(t))$ is
uniformly bounded on $[0, T]$ (see \cite{S}). We only consider the
singularity of the mean curvature flow, we assume that it blows up
at $ T$. The monotonicity formula for the mean curvature flow was
essentially proved by Huisken \cite{H2} and Hamilton \cite{Ha2}.
The weighted monotonicity formula for the symplectic mean
curvature flow was proved by Chen-Li in \cite{CL1}. In this
section, we prove one weighted monotonicity formula for the flow
(\ref{main}). Since the K\"ahler-Ricci flow exists $[0, T]$, thus
the injective of $(M, \overline{g}(t))$ is uniformly bounded from
below on $[0, T)$ \cite{TZ}. We can therefore choose a cut-off
function on $(M, \overline{g}(t))$ to prove the monotonicity
formula for the symplectic K\"ahler-Ricci mean curvature flow.

Let $H(X, X_0, t, t_0)$ be the backward heat kernel on
${\mathbb{R}}^4$. Define
$$\rho(X,  t)=(4\pi (t_0-t))H(X, X_0, t, t_0)=\frac{1}{ 4\pi
(t_0-t)} \exp{-\frac{|X-X_0|^2}{4(t_0-t)}}
$$ for $t<t_0$, such that
$$\frac{d}{dt}\rho=-\Delta\rho-\rho\left(\left|H+\frac{(X-X_0)^\perp}{2(t_0-t)}
\right|^2-|H|^2\right).
$$ where $(X-X_0)^\perp$ is the normal component of $X-X_0$.

Let $i_M$ be the lower bound of the injective radius of $(M,
\overline{g}(t))$. We choose a cut-off function $\phi\in
C^\infty_0(B_{2\tilde{r}}(X_0))$ with $\phi\equiv 1$ in
$B_{\tilde{r}}(X_0)$, where $X_0\in M$, $0<2\tilde{r}<i_M$. Choose
a normal coordinates in $B_{2\tilde{r}}(X_0)$ in $(M,
\overline{g}(t))$ and express $F$ using the coordinates $(F^1,
F^2, F^3, F^4)$ as a function in $\bf{R}^4$. Set $v(x,
t)=e^{\bar{R}_0t}\cos\alpha(x, t)$, where $\overline{R}_0=\max\{0,
-\overline{R}\}$ and $\overline{R}$ is the scalar curvature of
$M$. We define
$$\Phi(F, X_0, t, t_0)=\int_{\Sigma_t}\frac{1}{v}\phi(F)\rho(F, t)d\mu_t.
$$

\begin{proposition}\label{p1}
There are positive constants $c_1$ and $c_2$  depending only on
$(M,
\overline{g}(t))$, $F_0$ and $\tilde{r}$ such that
\begin{eqnarray}\label{e10}
\frac{\partial}{\partial t} e^{c_1\sqrt{t_0-t}}\Phi(F, X_0, t,
t_0)&\leq&
-e^{c_1\sqrt{t_0-t}}\left(\int_{\Sigma_t}\frac{1}{v}\phi\rho(F,
t)|H+\frac{(F-X_0)^\perp}{2(t_0-t)}|^2
d\mu_t\right.\nonumber\\&&\left.+\int_{\Sigma_t}\frac{1}{v}|\overline\nabla
J_{\Sigma_t}|^2\phi\rho(F, t)
d\mu_t+\int_{\Sigma_t}\frac{2|\nabla\cos\alpha|^2}{v^3}\phi\rho(F,t
) d\mu_t\right) \nonumber\\&&+c_2(t_0-t).
\end{eqnarray}
\end{proposition}

{\it Proof.} By the evolution equation of the K\"ahler angle in Theorem \ref{th2}, we have
$$(\frac{\partial}{\partial t}-\Delta)\frac{1}{v}\leq -\frac{|\overline\nabla J|^2}{v}
-\frac{2|\nabla v|^2}{v^3}.$$ So,
\begin{eqnarray*}
\frac{\partial}{\partial t}\Phi(F, X_0, t,
t_0)&=&\int_{\Sigma_t}\frac{\partial}{\partial
t}\frac{1}{v}\phi\rho(F,
t)d\mu_t+\int_{\Sigma_t}\frac{1}{v}\frac{\partial}{\partial
t}\rho(F, t)\phi
d\mu_t\\&&-\int_{\Sigma_t}\frac{1}{v}\phi\rho(F, t)|H|^2
d\mu_t-\int_{\Sigma_t}\frac{1}{v}\phi(\tilde{R}-r)\rho(F,
t)\\&\leq&\int_{\Sigma_t}\frac{1}{v}\phi(\frac{\partial}{\partial
t}+\Delta)\rho(F,
t)d\mu_t-\int_{\Sigma_t}\frac{1}{v}\phi\rho(F, t)|H|^2 d\mu_t
\\&&-\int_{\Sigma_t} \frac{|\overline\nabla J_{\Sigma_t}|^2}{v}\phi\rho(F, t)
d\mu_t-\int_{\Sigma_t}\frac{2|\nabla v|^2}{v^3}
\phi\rho(F, t)d\mu_t
\\
&&+\int_{\Sigma_t}\frac{1}{v}\Delta\phi\rho(F,
t)d\mu_t+2\int_{\Sigma_t}\frac{1}{v}\nabla\phi\nabla\rho(F,
t)d\mu_t\\&&-\int_{\Sigma_t}\frac{1}{v}(\tilde{R}-r)\phi\rho(F,
t) d\mu_t.
\end{eqnarray*}

In $[0, T)$, $\tilde{R}$ is uniformly bounded, thus
$$\left|\int_{\Sigma_t}\frac{1}{v}\phi(\tilde{R}-r)\rho(F, t)\right|\leq C$$
and by (\ref{e100})
$$\frac{\partial}{\partial t}\int_{\Sigma_t}d\mu_t\leq C\int_{\Sigma_t} d\mu_t,$$ so,
$$Area(\Sigma_t)\leq e^{CT} Area(\Sigma_0).$$
where the constant $C$ depends on $(M, \overline{g}(t))$.

Straight computation leads to
\begin{eqnarray*}
\frac{\partial}{\partial t}\rho(F,
t)=(\frac{1}{t_0-t}-\frac{\langle H,
F-X_0\rangle}{2(t_0-t)}-\frac{|F-X_0|^2}{4(t_0-t)^2})\rho(F, t)
\end{eqnarray*}
and
\begin{eqnarray*}
&&\Delta\exp(-\frac{|F-X_0|^2}{4(t_0-t)})\\&=&
\exp(-\frac{|F-X_0|^2}{4(t_0-t)})(\frac{|\langle F-X_0, \nabla
F\rangle|^2}{4(t_0-t)}-\frac{\langle F-X_0, \Delta
F\rangle}{2(t_0-t)}-\frac{|\nabla F|^2}{2(t_0-t)}).
\end{eqnarray*}
Notice that, in the induced metric on $\Sigma_t$,
$$|\nabla F|^2=2~~~~~~~~~~~~~~~~~~~~~~~~~~~~~and~~~~~~~~~~~~~~~~~~~~~~~~~~~~~~~~~~~~\Delta
F^\alpha=H^\alpha-g^{ij}\overline{\Gamma}^\alpha_{\rho\sigma}\frac{\partial
F^\rho}{\partial x^i}\frac{\partial F^\sigma}{\partial x^j},
$$ then we have
\begin{eqnarray*}
&&(\frac{\partial}{\partial t}+\Delta)\rho(F,
t)\\&=&-(\frac{\langle F-X_0,
H\rangle}{(t_0-t)}+\frac{|(F-X_0)^\bot|^2}{4(t_0-t)^2}+\frac{\langle
F-X_0, g^{ij}\overline{\Gamma}^\alpha_{\rho\sigma}\frac{\partial
F^\rho}{\partial x^i}\frac{\partial F^\sigma}{\partial
x^j}e_\alpha\rangle}{2(t_0-t)})\rho(F, t).
\end{eqnarray*}

Note that $\Delta\phi=0, \nabla\phi=0$ in $B_r(X_0)$, we can see
that
$$|\Delta\phi\rho(F, t)|\leq C~~~~~~~~~and~~~~~~~~~~|\nabla\phi\nabla\rho(F, t)|\leq
C.$$ Hence
$$\int_{\Sigma_t}\frac{1}{v}\Delta\phi\rho(F, t)d\mu_t\leq C\int_{\Sigma_t}d\mu_t\leq C$$
$$\int_{\Sigma_t}\frac{1}{v}\nabla\phi\nabla\rho(F, t)d\mu_t\leq C\int_{\Sigma_t}d\mu_t\leq C.$$

Since we choose a normal coordinates in $B_{2\tilde {r}}(X_0)$ in
$(M, \overline{g}(t))$, we have
$\overline{\Gamma}^\alpha_{\rho\sigma}(X_0, t)=0$, and
$|g^{ij}\overline{\Gamma}^\alpha_{\rho\sigma}\frac{\partial
F^\rho}{\partial x^i}\frac{\partial F^\sigma}{\partial x^j}|\leq
C|F-X_0|$, thus
\begin{eqnarray*}
\frac{\langle F-X_0,
g^{ij}\overline{\Gamma}^\alpha_{\rho\sigma}\frac{\partial
F^\rho}{\partial x^i}\frac{\partial F^\sigma}{\partial
x^j}e_\alpha\rangle}{2(t_0-t)}\leq C\frac{|F-X_0|^2}{2(t_0-t)}.
\end{eqnarray*} Hence
\begin{eqnarray*}
\frac{\langle F-X_0,
g^{ij}\overline{\Gamma}^\alpha_{\rho\sigma}\frac{\partial
F^\rho}{\partial x^i}\frac{\partial F^\sigma}{\partial
x^j}e_\alpha\rangle}{2(t_0-t)}\rho(F, t)\leq c_2\frac{\rho(F,
t)}{\sqrt{t_0-t}}+c_3.
\end{eqnarray*}

It concludes that \begin{eqnarray*} \frac{\partial}{\partial
t}\Phi(F, X_0, t, t_0)&\leq&
-\int_{\Sigma_t}\frac{1}{v}\phi\rho(F,
t)|H+\frac{(F-X_0)^\bot}{2(t_0-t)}|^2 d\mu_t\\&& -\int_{\Sigma_t}
\frac{|\overline\nabla J_{\Sigma_t}|^2}{v}\phi\rho(F, t)
d\mu_t-\int_{\Sigma_t}\frac{2|\nabla v|^2}{v^3}
\phi\rho(F, t)d\mu_t\\
&&+\frac{c_1}{\sqrt{t_0-t}}\Phi(F, X_0, t, t_0)+c_2.
\end{eqnarray*} The proposition follows. \hfill Q. E. D.

\section{ no type I singularity}

We assume in this section that the K\"ahler-Ricci flow exists
globally. We study the singularity of the mean curvature flow if
it blows up at $T$. Since the K\"ahler-Ricci flow exists for all
time, the Ricci curvature is uniformly bounded, and the evolution
equation of the K\"ahler angle along the K\"ahler-Ricci mean
curvature flow is the same as that of the mean curvature flow in
K\"ahler-Einstein surface (c.f (\ref{cos}) and Proposition 3.2 in
\cite{CL1}), thus the analysis of the singularity of the mean
curvature flow is the same as that of the mean curvature flow in
K\"aher-Einstein surface \cite{CL1}. For completeness, we give
some details below.

We recall the classification of the singularities of the mean
curvature flow. We say the mean curvature flow has type I
singularity at $T>0$, if
$$\overline{\lim_{t\to T}}(T-t)\max_{\Sigma_t}|A|^2\leq C,$$ for
some positive constant $C$. Otherwise, we say the mean curvature
flow has type II singularity.

\begin{lemma}\label{l1}
Let $U(t)=\max_{\Sigma_t}|A|^2$. If the mean curvature flow blows up at $T>0$,
then there is a positive constant $c$ depending only on the bound of the curvature $(M, \overline{g}(t))$ such that,
if $0<T-t<\pi/16\sqrt{c}$, the function $U(t)$ satisfies
$$U(t)\geq\frac{1}{4\sqrt{2}(T-t)}.$$
\end{lemma}

{\it Proof.} By lemma\ref{A} and the parabolic maximum principle, we have
\begin{eqnarray*}
\frac{\partial}{\partial t} U(t) &\leq& 2(U(t))^2+c_1U(t)+c_2\sqrt{U(t)}\\
&\leq& 4(U(t))^2+4c,
\end{eqnarray*} where $c_1, c_2$ are constants which depend only on the bounds of the curvature and its covariant derivatives of $(M, \overline{g}(t))$.
This implies the desired inequality. \hfill Q. E. D.

\begin{theorem}\label{th7}  Assume that the K\"ahler-Ricci flow exists globally.
The symplectic K\"ahler-Ricci mean curvature flow has no type I singularity at any $T>0$.
\end{theorem}
{\it Proof.}  Suppose that the mean curvature flow has a type I singularity at $t_0>0$. Assume that
$$\lambda_k^2=|A|^2(x_k, t_k)=\max_{t\leq t_k}|A|^2$$ and $x_k\to p\in \Sigma$, $t_k\to t_0$ as $k\to\infty$.
We choose a local coordinate system on $(M, \overline{g}((t))$
around $F(p, t_0)$ such that $F(p, t_0)=0$. And we rescale the
mean curvature flow,
$$F_k(x, t)=\lambda_k(F(x, \lambda_k^{-2}t+t_k)-F(p, t_k)),~~~~~~~~~~~~t\in [-\lambda_k^2t_k, 0].$$ Denote by $\Sigma^k_t$
the rescaled surface $F_k(\cdot, t)$.  By Lemma \ref{l1}, we have
$$\frac{C}{t_0-t_k}\geq |A|^2(x_k, t_k)\geq\frac{c}{t_0-t_k}$$ for some uniform constants $c$ and $C$ independent of $k$. Therefore,
$$|A_k|^2(F(x_k,t_k)-F(p,t_k))=|A|^2(F_k(x_k,0))=\lambda_k^{-2}|A|^2(x_k, t_k)\geq \lambda_k^{-2}\frac{c}{(t_0-t_k)}.$$ Since the mean curvature flow has type
I singularity at $t_0>0$, we have
$$\lambda_k^{2}(t_0-t_k)\leq C.$$ So,
$$|A_k|^2(0)\geq c>0,$$ for some uniform constants $c$. It is easy to see that
$$|A_k|^2(x, t)\leq 1,$$ thus there is a subsequence of $F_k$ which we also denote by $F_k$, such that
$F_k\to F_\infty$ in any ball $B_R(0)\subset {\textbf{R}}^4$, and $F_\infty$ satisfies
\begin{equation}\label{e14}\frac{\partial F_\infty}{\partial t}=H_\infty~~~~~~~~~~{\rm with}~~~~~~~~~~~~~1\geq |A_\infty|(0)\geq c>0.\end{equation}
By the monotonicity
formula (\ref{e10}), we know that $\lim_{t\to t_0} e^{c_1\sqrt{t_0-t}}\Phi$ exists. Using the equality
\begin{eqnarray*}
&&\int_{\Sigma^k_t}\frac{1}{v}\frac{1}{0-t}\phi_{\lambda_kr}(F_k)\exp(-\frac{|F_k+\lambda_k F(p, t_k)|^2}{4(0-t)})d\mu_t^k \\
&&=\int_{\Sigma_{t_k+\lambda_k^{-2}t}}\frac{1}{v}\phi(F)\frac{1}{t_k-(t_k+\lambda_k^{-2}t)}
\exp(-\frac{|F(x, t_k+\lambda_k^{-2}t)|^2}{4(t_k-(t_k+\lambda_k^{-2}t))})d\mu_t,
\end{eqnarray*} we can get that, for any $-\infty<s_1<s_2<0$,
\begin{eqnarray*}
&&e^{c_1\sqrt{t_k-(t_k+\lambda_k^{-2}s_2)}}\int_{\Sigma^k_{s_2}}
 \frac{1}{v}\frac{1}{0-s_2}\phi_{\lambda_kr}(F_k)\exp(-\frac{|F_k+\lambda_k F(p, t_k)|^2}{4(0-s_2)})d\mu_{s_2}^k\\&&-
 e^{c_1\sqrt{t_k-(t_k+\lambda_k^{-2}s_1)}}\int_{\Sigma^k_{s_1}}
 \frac{1}{v}\frac{1}{0-s_1}\phi_{\lambda_kr}(F_k)\exp(-\frac{|F_k+\lambda_k F(p, t_k)|^2}{4(0-s_1)})d\mu_{s_1}^k \\&&\rightarrow 0
 ~~~~~~~~~~~~{\rm as}~~~~~~~~~~k\to\infty.
\end{eqnarray*} Integrating (\ref{e10}) from $t_k+\lambda_k^{-2}s_1$ to $t_k+\lambda_k^{-2}s_2$, we get that
\begin{eqnarray}\label{e12}
|\overline{\nabla} J_{\Sigma_{t_k}}|^2\longrightarrow 0,
\end{eqnarray} and
\begin{eqnarray}\label{e13}
|H_k+\frac{(F_k+\lambda_kF(p, t_k))^\bot}{2(0-t)}|^2\longrightarrow 0~~~~~~~~~{\rm as}~~~~~~~~~~~k\to\infty.
\end{eqnarray} By (\ref{e12}), we get that
$$DJ_\infty\equiv 0,$$ where $D$ is the derivative in ${\bf{R}}^4$, which implies that
$$H_\infty\equiv 0.$$ Since
\begin{eqnarray*}
|F(p, t_k)|\leq\int_{t_k}^{t_0}|\frac{\partial F_k}{\partial t}|dt\leq \int_{t_k}^{t_0}|H|dt\leq C\sqrt{t_0-t_k}\leq\frac{C}{\lambda_k},
\end{eqnarray*} thus $\lambda_k F(p, t_k)\to q$ as $k\to\infty.$
Hence by (\ref{e13}), we have
$$(F_\infty+q)^\perp\equiv 0,$$ this implies that, for $\alpha=3, 4$,
$$\det((h_\infty)^\alpha_{ij})=0.$$ Since $H_\infty=0$, we also have, for $\alpha=3, 4$,
$$tr((h_\infty)^\alpha_{ij})=0.$$ This yields that $(h_\infty)^\alpha_{ij}=0$ for all $i,j=1, 2, \alpha=3, 4$ which  implies that
$|A_\infty|\equiv 0$. This contradicts with (\ref{e14}). \hfill Q. E. D.

\section{graph case}

In this section we study the K\"ahler-Ricci mean curvature flow
(\ref{main}) in a special case. Suppose that  $M$ is a product of
compact Riemann surfaces $M_1, M_2$, i.e, $(M,
\overline{g})=(M_1\times M_2,
\overline{g}_1\oplus\overline{g}_2)$. We denote by $r_1, r_2$ the
average scalar curvature of  $M_1, M_2$, we assume that $r_1=r_2$.
Then the K\"ahler-Ricci flow on $M$ can be split into the
K\"ahler-Ricci flow on $M_1$, $M_2$ respectively. It is well known
that the K\"ahler-Ricci flow on surface exists for long time and
converges to the surface with constant curvature at infinity.
Suppose that $\Sigma$ is a graph in $M=M_1\times M_2$. Recall the
definition of the graph in \cite{CLT}.  A surface $\Sigma$ is a
graph in $M_1\times M_2$ if $v=\langle e_1\wedge e_2,
\omega_1\rangle\geq c_0>0$ where $\omega_1$ is a unit K\"ahler
form on $M_1$, and $\{e_1, e_2\}$ is an orthonormal frame on
$\Sigma$. In this section, we used some ideas in \cite{CLT} and
\cite{Wa}.

\begin{theorem}\label{th9}
Let $(M_1, \overline{g}_1, \omega_1)$ and $(M_2, \overline{g}_2, \omega_2)$ be
Riemann surfaces which have the same average scalar curvature. Suppose that $M_1\times
M_2$ evolves along the Ricci flow flow and $\Sigma_0$ evolves along the
mean curvature in $M_1\times M_2$. If $v(\cdot,
0)>\frac{\sqrt{2}}{2}$, then the K\"ahler-Ricci mean curvature flow exists for
all time.
\end{theorem}

{\it Proof.} Because $M_1$ and $M_2$ have the same average scalar
curvature, the metric
$\overline{g}=\overline{g}_1\oplus\overline{g}_2$ on $M$ evolves
along the K\"ahler-Ricci flow is equivalent to $\overline{g}_1$
and $\overline{g}_2$ evolves along the Ricci flow respectively.
When $n=2$, $\overline{R}_{ij}=\frac{1}{2}\overline{R}g_{ij}$.
Thus, $\overline{g}_1(t)$,  $\overline{g}_2(t)$ satisfy the
evolution equations:
\begin{eqnarray}\label{e22}
\left\{ \begin{array}{clcr} \frac{\partial }{\partial t}(\overline{g}_1)_{ij} &= &-\frac{1}{2}(\overline{R}_1-r)(\overline{g}_1)_{ij}\nonumber\\
 (\overline {g}_1)_{ij}(0) &= &(\overline{g}_1)_0
\end{array}\right.,
\end{eqnarray}
and
\begin{eqnarray}\label{e23}
\left\{\begin{array}{clcr}\frac{\partial}{\partial t}(\overline{g}_2)_{ij} &= &-\frac{1}{2}(\overline{R}_2-r)(\overline{g}_2)_{ij}\nonumber\\
(\overline {g}_2)_{ij}(0) &= &(\overline{g}_2)_0
\end{array}\right..
\end{eqnarray} By the work of Hamilton \cite{Ha3} and Chow \cite{Ch}, we know that,
 for any initial metrics the flows exist for long time and converge to $M^\infty_1\times M^\infty_2$ with constant curvatures at infinity.

Choose an orthonormal basis $\{e_1, e_2, e_3, e_4\}$ on $M$ along
$\Sigma_t$ such that $\{e_1, e_2\}$ is the basis of $\Sigma_t$.
Set $u_1=\langle e_1\wedge e_2, \omega_1+\omega_2\rangle$ and
$u_2=\langle e_1\wedge e_2, \omega_1-\omega_2\rangle$ where
$\omega_2$ is an unit K\"ahler form on $M_2$. Since both
$\omega_1+\omega_2$ and $\omega_1-\omega_2$ are parallel K\"ahler
forms on $M_1\times M_2$, we see that Theorem \ref{th2} is
applicable. Therefore,
\begin{eqnarray}\label{e20}
(\frac{d}{dt}-\Delta)u_1=[(h^3_{2k}+h^4_{1k})^2+(h^4_{2k}-h^3_{1k})^2]u_1+\overline{R}_1u_1(1-u_1^2).
\end{eqnarray} By switching $e_3$ and $e_4$, we get that
\begin{eqnarray}\label{e21}
(\frac{d}{dt}-\Delta)u_2=[(h^4_{2k}+h^3_{1k})^2+(h^3_{2k}-h^4_{1k})^2]u_2+\overline{R}_2u_2(1-u_2^2).
\end{eqnarray} It is easy to see that
$\langle e_1\wedge e_2, \omega_1\rangle^2+\langle e_1\wedge e_2,
\omega_2\rangle^2\leq 1$. The initial condition $v(x, 0)\geq
\frac{\sqrt{2}}{2}$ implies that $u_i(x, 0)\geq v(x,
0)-\frac{\sqrt{2}}{2}\geq c_0>0, i=1, 2$. By (\ref{e20}) and
(\ref{e21}), applying the maximum principle for parabolic
equations,we see that $u_i(x, t)$ have positive lower bounds at
any finite time. Suppose that $u_i\geq \delta$ for $0\leq t<t_0$.
Then we claim that the flow $F$ can be extended smoothly to
$t_0+\varepsilon$ for some $\varepsilon$.

Set $u=u_1+u_2$. Adding (\ref{e20}) into (\ref{e21}), we get that
\begin{eqnarray}\label{e26}
(\frac{d}{dt}-\Delta)u&=&u|A|^2+2(u_1-u_2)h^3_{2k}h^4_{1k}-2(u_1-u_2)h^3_{1k}h^4_{2k}
\nonumber\\&&+ \overline{R}_1u_1(1-u_1^2)+\overline{R}_2u_2(1-u_2^2).
\end{eqnarray} Since $u\geq 2\delta+|u_1-u_2|$, using the Cauchy-Schwarz inequality, we get that
\begin{equation}\label{e22}
(\frac{d}{dt}-\Delta)u\geq 2\delta|A|^2-C,
\end{equation} where $C$ is the lower bound of the scalar curvature of $(M, {\overline{g}}_1(t)\oplus{\overline{g}}_2(t))$.

Assume that $(X_0,
t_0)$ is a singularity point. As in the proof of Proposition \ref{p1}, we can derive a weighted monotonicity formula for
$\int_{\Sigma_t}\phi\frac{1}{u}\rho(F, X_0, t, t_0) d\mu_t $, where $\phi$ is the cut-off function in Proposition \ref{p1}.
\allowdisplaybreaks \begin{eqnarray}\label{e3} \lefteqn{
\frac{d}{dt}\int_{\Sigma_t}\phi\frac{1}{u}\rho(F, X_0, t, t_0)d\mu_t }\nonumber\\
&\leq&\int_{\Sigma_t}\phi\rho\Delta \frac{1}{u}d\mu_t
-2\delta\int_{\Sigma_t}\phi\frac{|A|^2}{u^2}\rho d\mu_t
-2\int_{\Sigma_t}\phi\frac{|\nabla u|^2}{u^3}\rho d\mu_t
\nonumber\\&&+\int_{\Sigma_t}\phi\frac{C}{u^2}d\mu_t-\int_{\Sigma_t}\phi\frac{1}{u}\left(\Delta\rho
+\left(\left|H+\frac{(F-X_0)^\perp}{2(t_0-t)}\right|^2-|H|^2
\right)\rho\right)d\mu_t\nonumber\\
&&-\int_{\Sigma_t}\phi\frac{|H|^2}{u}\rho d\mu_t-\int_{\Sigma_t}\phi\frac{\tilde{R}-r}{u}\rho d\mu_t
\nonumber\\
&\leq&-\int_{\Sigma_t}\phi\rho\left(\frac{2}{u^3}|\nabla u|^2
+\frac{1}{u}\left|H+\frac{(F-X_0)^\perp}{2(t_0-t)}\right|^2
+2\delta\frac{|A|^2}{u^2} d\mu_t\right)\nonumber\\
&&-\int_{\Sigma_t}\phi\frac{\tilde{R}-r}{u}\rho d\mu_t+\int_{\Sigma_t}
\Delta\phi\frac{1}{u}\rho d\mu_t+2\int_{\Sigma_t}\frac{1}{u}\nabla\phi\cdot\nabla\rho d\mu_t
\nonumber\\&\leq& C- 2\delta\int_{\Sigma_t}\phi\frac{|A|^2}{u^2}\rho (F,
X_0, t, t_0)d\mu_t,
\end{eqnarray} where $C$ depends on the scalar curvature of $\overline{R}_1, \overline{R}_2$ and the bound of $|\nabla\phi|, |\Delta\phi|$. From this we see that $\lim_{t\to
t_0}\int_{\Sigma_t}\phi\frac{1}{u}\rho d\mu_t$ exists.

Let $0<\lambda_i\to\infty$ and let $F_i$ be the blow up sequence:
$$F_i(x, s)=\lambda_i(F(x, t_0+\lambda_i^{-2}s)-X_0).$$
Let $d\mu^i_s$ denote the induced volume form on $\Sigma^i_s$ by
$F_i$. It is obvious that,
$$\int_{\Sigma_t}\phi\frac{1}{u}\rho(F, X_0, t, t_0) d\mu_t
=\int_{\Sigma^i_s}\phi\frac{1}{u}\rho(F_i, 0, s, 0) d\mu^i_s.
$$ Therefore we get that, \allowdisplaybreaks
\begin{eqnarray*}
\lefteqn{
\frac{d}{ds}\int_{\Sigma^i_s}\phi\frac{1}{u}\rho(F_i, 0, s,
0)d\mu^i_s}\\ &\leq& \frac{C}{\lambda_i^2}-2\delta\int_{\Sigma^i_s}\phi\frac{|A_i|^2}{u^2}\rho
(F_i, 0, s, 0)d\mu_s^i
\end{eqnarray*} Note that $t_0+\lambda_i^{-2}s\to t_0$ for any
fixed $s$ as $i\to\infty$ and that $\lim_{t\to
t_0}\int_{\Sigma_t}\phi\frac{1}{u}\rho d\mu_t$ exists. By the
above monotonicity formula, we have, for any fixed $s_1$
and $s_2$,
\begin{eqnarray*}
0&\leftarrow&\int_{\Sigma^i_{s_1}}\phi\frac{1}{u}\rho(F_i, 0,
s_1,
0)d\mu^i_{s_1}-\int_{\Sigma^i_{s_2}}\phi\frac{1}{u}\rho(F_i,
0, s_2, 0) d\mu^i_{s_2}\\ &\geq&
2\delta\int_{s_1}^{s_2}
\int_{\Sigma^i_s}\phi\frac{|A_i|^2}{u^2}\rho (F_i, 0, s,
0)d\mu_s^i.
\end{eqnarray*} Since $u$ is bounded below, we have
$$\int_{s_1}^{s_2}\int_{\Sigma_s^i}|A_i|^2\rho(F_i, 0, s, 0)\to 0
~~~~~~~~{\rm as}~~~~~~~~i\to\infty.
$$ Therefore, for any ball $B_R(0)\subset R^{4}$,
\begin{equation}\label{e24} \int_{\Sigma_{s_i}^i\cap B_R(0)}|A_i|^2\to 0 ~~~~~~~~{\rm
as}~~~~~~~i\to \infty.
\end{equation}

Because
$u$ has a positive lower bound, we see that
$\Sigma_t$ can locally be written as the graph of a map $f_t:\Omega\subset M_1\to
M_2$ with uniformly bounded $|df_t|$. Consider the blow up of
$f_{t_0+\frac{s_i}{\lambda_i^2}}$ ,
$$f_i(y)=\lambda_i f_{t_0+\lambda_i^{-2}s_i}(\lambda_i^{-1}y).
$$ It is clear that $|d f_i|$ is also uniformly bounded
and $\lim_{i\to\infty} f_i(0)=0$. By Arzzela theorem, $ f_i\to
f_\infty$ in $C^\alpha$ on any compact set. By the inequality
$(29)$ in \cite{I}, we have
$$ |A_i|\leq |\nabla df_i|\leq C(1+|df_i|^3)|A_i|,
$$ where $\nabla df_i$ is measured with respect to the induced
metric on $\Sigma^i_{s_i}$. From equation (\ref{e24}) it follows that, for any ball $B_R(0)\subset R^{4}$,
$$ \int_{\Sigma_{s_i}^i\cap B_R(0)}|\nabla df_i|^2\to 0
~~~~~~~~{\rm as}~~~~~~~~i\to\infty,
$$ which implies that $f_i\to f_\infty$ in $C^\alpha\cap W^{1,
2}_{loc}$ and the second derivative of $f_\infty$ is $0$. It is then clear that
$\Sigma^i_{s_i}\to\Sigma^{\infty}$ and  $\Sigma^\infty$ is the graph of a
linear function. Therefore,
$$\lim_{i\to\infty}\int\phi\rho(F_i, 0, s_i, 0)d\mu^i_{s_i}=\int\rho(F_\infty, 0, -1,
0)d\mu^\infty=1,
$$ We therefore have
\begin{eqnarray}\label{e25}\lim_{t\to t_0}\int\rho(F, X_0, t,
t_0)&=&\lim_{i\to\infty}\int\phi\rho(F, X_0,
t_0+\lambda_i^{-2}s_i, t_0)\nonumber \\
&=&\lim_{i\to\infty}\int\phi\rho(F_i, 0, s_i, 0)d\mu^i_{s_i}=1.
\end{eqnarray} By the White's regularity theorem \cite{W}, we know that $(X_0, t_0)$ is a
regular point. This proves the theorem. \hfill Q. E. D.

\begin{theorem}\label{th10}
Under the same assumption as in Theorem \ref{th9}. If the scalar
curvature of $M_1, M_2$ is positive, then it converges to a
totally geodesic surface at infinity.
\end{theorem}

{\it Proof.} Since the scalar curvature of $M$ is positive, by (\ref{e20}) and
(\ref{e21}),
\begin{eqnarray*}
(\frac{\partial}{\partial t}-\Delta)(1-u_1)\leq-\overline{R}_1u_1(1+u_1)(1-u_1)\leq -c_1\overline{R}_1(1-u_1),
\end{eqnarray*} where $c_1$ depends only on the lower bound of $u_1$. Applying the maximum principle, we get that $1-u_1\leq ce^{-c_1\bar{R}_1t}$.  Similarly,
$1-u_2\leq ce^{-c_2\bar{R}_2t}$, where $c_1,c_2$ depends only on
the lower bound of $u_1,u_2$.  By (\ref{e26}), for any
$\varepsilon>0$, there exists $T$ such that as $t>T$,
$u_1>1-\varepsilon$, $u_2>1-\varepsilon$, $|u_1-u_2|<\varepsilon$
and
\begin{eqnarray*}\label{e27}
(\frac{\partial}{\partial t}-\Delta)u\geq (1-\varepsilon)|A|^2.
\end{eqnarray*}
From (\ref{e6}) we see that
\begin{eqnarray*}
(\frac{\partial}{\partial t}-\Delta)|A|^2 \leq -2|\nabla A|^2+C_1|A|^4+C_2|A|^2+C_3,
\end{eqnarray*} where $C_1, C_2, C_3$ are constants that depend on the bounds of the curvature
tensor and its covariant derivatives of $(M,\overline{g}(t))$.

Let $p>1$ be a constant to be fixed later. Now we consider the function $\frac{|A|^2}{e^{pu}}$.
\begin{eqnarray*}
(\frac{\partial}{\partial t}-\Delta)\frac{|A|^2}{e^{pu}}&=& 2\nabla(\frac{|A|^2}{e^{pu}})\cdot\frac{\nabla e^{pu}}{e^{pu}}\\
&&+ \frac{1}{e^{2pu}}[e^{pu}(\frac{\partial}{\partial t}-\Delta)|A|^2-|A|^2(\frac{\partial}{\partial t}-\Delta)e^{pu}]\\ &\leq&
2p\nabla(\frac{|A|^2}{e^{pu}})\cdot\nabla u \\&&+\frac{1}{e^{2pu}}[e^{pu}(C_1|A|^4+C_2|A|^2+C_3)-p|A|^2e^{pu}[(1-\varepsilon)|A|^2-p|\nabla u|^2]]
.
\end{eqnarray*}

From (\ref{e28}) and (\ref{e29}), it follows that,
\begin{eqnarray*}
|\nabla u_1|^2 &\leq& 2(1-u_1^2)((h^4_{1k})^2+(h^3_{2k})^2),\\
|\nabla u_2|^2 &\leq& 2(1-u_2^2)((h^3_{1k})^2+(h^4_{2k})^2).
\end{eqnarray*} So, for $t$ is sufficiently large, we have
\begin{eqnarray*}
|\nabla u_1|^2\leq \varepsilon |A|^2,~~~~~~~~~~~~~~~~~~~~~~~~~~~~~~~~~~~~~~~~|\nabla u_2|^2\leq\varepsilon |A|^2.
\end{eqnarray*}
and $$|\nabla u|^2\leq 2(|\nabla u_1|^2+|\nabla u_2|^2)\leq 4\varepsilon |A|^2.$$
Therefore,
\begin{eqnarray*}
(\frac{\partial}{\partial t}-\Delta)\frac{|A|^2}{e^{pu}} &\leq& 2p\nabla(\frac{|A|^2}{e^{pu}})\cdot\nabla u \\&&+
\frac{1}{e^{pu}}[(C_1-p(1-\varepsilon)+4p^2\varepsilon)|A|^4+C_2|A|^2+C_3].
\end{eqnarray*}
Set $p^2=1/\varepsilon$, then
\begin{eqnarray*}
C_1-p(1-\varepsilon)+4p^2\varepsilon &=&C_1-\varepsilon^{-\frac{1}{2}}+\varepsilon^{\frac{1}{2}}+4.
\end{eqnarray*} As $t$ is sufficiently large, i.e. $\varepsilon$ is sufficiently close to $0$, we have
\begin{eqnarray*}
(C_1-\varepsilon^{-\frac{1}{2}}+\varepsilon^{\frac{1}{2}}+4)\leq -1.
\end{eqnarray*} So,
\begin{eqnarray*}
(\frac{\partial}{\partial t}-\Delta)\frac{|A|^2}{e^{pu}} &\leq& 2p\nabla(\frac{|A|^2}{e^{pu}})\cdot\nabla u
-\frac{|A|^4}{e^{pu}}+C_2\frac{|A|^2}{e^{pu}}+\frac{C_3}{e^{pu}}\\ &\leq& 2p\nabla(\frac{|A|^2}{e^{pu}})\cdot\nabla u
-\frac{|A|^4}{e^{2pu}}+C_2\frac{|A|^2}{e^{pu}}+\frac{C_3}{e^{pu}}
\end{eqnarray*}

Applying the maximum principle for parabolic equations, we conclude that $\frac{|A|^2}{e^{pu}}$ is uniformly bounded, thus
$|A|^2$ is also uniformly bounded. Thus $F(\cdot, t)$ converges to $F_\infty$ in $C^2$ as $t\to\infty$. Since
\begin{eqnarray*}
1-u_1\leq c e^{-c_1{\bar R}_1t},
\end{eqnarray*} and
\begin{eqnarray*}
1-u_2\leq c e^{-c_1{\bar R}_2t},
\end{eqnarray*}
So, we have $u_1\equiv 1$ and $u_2\equiv 2$ at infinity. By
(\ref{e26}), we see that the second fundamental form is zero
identically, in other words, the limiting surface $F_\infty$ is
totally geodesic.

\hfill Q. E. D.

\section{stability of K\"ahler-Ricci mean curvature flow}

Let $(M, J)$ be a K\"ahler surface with $c_1(M)>0$. Now we choose
$r=2$. Suppose $(M, J)$ is pre-stable and the Futaki invariant of
the class $2\pi c_1(M)$ vanishes. Chen-Li \cite{ChL} proved that
for any $\gamma$, there exists a small positive constant
$\varepsilon(\gamma)$ such that for any metric $\overline{g}$ in
the subspace of K\"ahler metrics
$$\{\omega_{\bar{g}}\in2\pi c_1(M)| |Rm|(\omega_{\bar{g}})\leq\gamma, |Ric(\omega_{\bar{g}})-\omega_{\bar{g}}|\leq\varepsilon\},
$$
the K\"ahler Ricci flow with the initial metric $\omega_{\bar{g}}$ will converge exponentially fast to a K\"ahler-Einstein metric.

That is
$$
\|\overline{Ric}(g(t, \cdot))-\overline{g}(t, \cdot)\|_{\bar{g}(t, \cdot)}\leq C\varepsilon e^{-\beta t},
$$
for some positive constants $C$ and $\beta$ which depend only on
$\gamma$ and $\varepsilon$. In this case, we show that, if the
initial surface is sufficiently close to a holomorphic curve, then
the K\"ahler-Ricci mean curvature flow exists for all time (the
idea is similar to that in \cite{HL}) and converges to $(M,
\overline{g}_\infty, \Sigma_\infty)$ at infinity.

Note that, since $\Sigma_0$ is smooth, it is well-known that
$$\lim_{r\rightarrow 0}\int_{\Sigma_0}\phi(F)\frac{1}{4\pi r^2}
e^{\frac{|F-X_0|^2}{4r^2}}d\mu_0=1
$$ for any $X_0\in\Sigma_0$. So we can find a sufficiently small
$r_0$ such that
\begin{equation}\label{e95}\int_{\Sigma_0}\phi(F)\frac{1}{4\pi
r_0^2}e^{\frac{|F-X_0|^2}{4r_0^2}}d\mu_0\leq 1+\varepsilon_0/2
\end{equation} for all $X_0\in M$, where $\varepsilon_0$ is the constant in White's Theorem.

\begin{theorem}\label{th8} There exist sufficiently small constant $\varepsilon_1, \varepsilon_2$ such that, if
$(\varepsilon_1+\varepsilon_2)/r_0^6\ll\varepsilon_0$ where $r_0$ is defined in (\ref{e95}) and $\varepsilon_0$ is a constant in White's theorem,
and at any time $t$,
 $$
\|\overline{Ric}(g(t, \cdot))-\overline{g}(t, \cdot)\|_{\bar{g}(t, \cdot)}\leq C\varepsilon_1 e^{-\beta t},
$$
and the K\"ahler angle of the initial surface satisfies
$\cos\alpha_0\geq 1-\varepsilon_2$,  then K\"ahler-Ricci mean
curvature flow exists globally and it converges to $(M,
\overline{g}_\infty, \Sigma_\infty)$ at infinity. Furthermore,
$\Sigma_\infty$ is the holomorphic curve in $(M,
\overline{g}_\infty)$ and $(M, \overline{g}_\infty)$ is
K\"ahler-Einstein surface.
\end{theorem}

{\it Proof.} Since along the Ricci flow, at any time $t$ we have,
\begin{eqnarray*}\label{e98}
|\overline{Ric}-\omega|\leq C\varepsilon_1 e^{-\beta t},
\end{eqnarray*} then at any time $t$,
\begin{eqnarray}
|\tilde{R}-2|\leq C\varepsilon_1 e^{-\beta t}.
\end{eqnarray}  By (\ref{cos}),  along the mean curvature flow,
$\cos\alpha$ increase, i.e, $\cos\alpha\geq 1-\varepsilon_2.$ Using the equation of (\ref{sin}) we can obtain that
\begin{eqnarray*}
(\frac{\partial}{\partial t}-\Delta)\sin^2\alpha/2\leq-2c\varepsilon_2(1+\varepsilon_2)\sin^2\alpha/2 .
\end{eqnarray*} So,
\begin{eqnarray}\label{e96}
\sin^2\alpha/2\leq \sin^2\alpha_0/2 e^{-ct}\leq \varepsilon_2 e^{-ct}.
\end{eqnarray}

By (\ref{e100}),
\begin{eqnarray*}
\frac{d}{dt}\int_{\Sigma_t}d\mu_t\leq\int_{\Sigma_t}|\tilde{R}-2|d\mu_t\leq C\varepsilon_1 e^{-\beta t} \int_{\Sigma_t}d\mu_t,
\end{eqnarray*} thus,
\begin{eqnarray*}
Area(\Sigma_t)\leq C Area(\Sigma_0).
\end{eqnarray*}

Because $\omega(t)$ is always closed, we can see that
$$
\int_{\Sigma_t}\cos\alpha d\mu_t = \int_{\Sigma_t}\omega
$$
is constant under the continuous deformation in $t$. Thus ,
\begin{eqnarray}\label{e16}
\frac{\partial}{\partial
t}\int_{\Sigma_t}(1-\cos\alpha)d\mu_t&=&-\int_{\Sigma_t} |H|^2 d\mu_t -\int_{\Sigma_t}
(\tilde{R}-2)d\mu_t.
\end{eqnarray}
Integrating the above inequality from $t$ to $t+1$ we obtain that
\begin{eqnarray*}
\int_t^{t+1}\int_{\Sigma_t} |H|^2 d\mu_t dt&\leq& \int_{\Sigma_t}\sin^2\alpha/2 d\mu_t+\int_t^{t+1}\int_{\Sigma_s}|\tilde{R}-2|d\mu_s ds \\
&\leq& C Area(\Sigma_0)(\varepsilon_2 e^{-ct} +\varepsilon_1 e^{-\beta t}) \\ &\leq& C(\varepsilon_1+\varepsilon_2) e^{-\lambda t},
\end{eqnarray*}  where $\lambda=\min\{c, \beta\}$. From this we can derive an $L^1$-estimate of the mean curvature vector.
\begin{eqnarray*}
\int_0^t\int_{\Sigma_s}|H|d\mu_s
ds&=&\sum_{k=0}^{t-1}\int_k^{k+1}\int_{\Sigma_s}|H|d\mu_s
ds\\&\leq&\sum_{k=0}^{t-1}(\int_k^{k+1}\int_{\Sigma_s}|H|^2 d\mu_s
ds)^{1/2}(\int_k^{k+1}Area\Sigma_t)^{1/2}\\&\leq&
C Area(\Sigma_0)^{1/2}\sum_{k=0}^{t-1}(\int_k^{k+1}\int_{\Sigma_s}|H|^2
d\mu_s ds)^{1/2}\\&\leq&
C(\varepsilon_1+\varepsilon_2)^{1/2}\sum_{k=0}^{t-1}e^{-\frac{\lambda k}{2}}\\&\leq&
C\frac{ (\varepsilon_1+\varepsilon_2)^{1/2}}{{1-e^{-\frac{\lambda}{2}}}},
\end{eqnarray*} where $C$ depends only on the area of the initial surface $\Sigma_0$.

Now we explore a possible singularity $(X_0, T)$. We study the density in White's local regularity theorem \cite{W}. Recall that
it is defined by
$$\Phi(X, X_0, t, t-r^2)=\int_{\Sigma_{t-r^2}}\phi(F)
\frac{1}{4\pi r^2}e^{-\frac{|F-X_0|^2}{4r^2}}d\mu_{t-r^2},$$ where $\phi$ is a cut off function around $X_0$ on $M$ such that $\phi\equiv 1$ in $B_r(X_0)$.
Differentiating this equation with respect to $t$ we get that
\begin{eqnarray}\label{e99}
\frac{\partial}{\partial
t}\int_{\Sigma_{t-r^2}}\phi(F)\frac{1}{4\pi
r^2}e^{-\frac{|F-X_0|^2}{4r^2}}d\mu_{t-r^2}&=&\int_{\Sigma_{t-r^2}}
\bigtriangledown\phi\cdot H \frac{1}{4\pi
r^2}e^{-\frac{|F-X_0|^2}{4r^2}}d\mu_{t-r^2}\nonumber\\&&-\int_{\Sigma_{t-r^2}}\frac{\phi}{8\pi
r^4}e^{-\frac{|F-X_0|^2}{4r^2}}\langle F-X_0, H\rangle
d\mu_{t-r^2}\nonumber\\&&-\int_{\Sigma_{t-r^2}}\frac{\phi}{4\pi
r^2}e^{-\frac{|F-X_0|^2}{4r^2}}|H|^2 d\mu_{t-r^2} \nonumber\\&&-\int_{\Sigma_{t-r^2}}\frac{\phi}{4\pi
r^2}e^{-\frac{|F-X_0|^2}{4r^2}} (\tilde{R}-2) d\mu_{t-r^2}.\nonumber\\
\end{eqnarray}
Integrating (\ref{e99}) from $r_0^2$ to $T$ we get that
\begin{eqnarray*}
&&\int_{\Sigma_{T-r_0^2}}\phi(F)\frac{1}{4\pi
r_0^2}e^{-\frac{|F-X_0|^2}{4r_0^2}}d\mu_{T-r_0^2}\leq\int_{\Sigma_0}\phi(F)\frac{1}{4\pi
r_0^2}e^{-\frac{|F-X_0|^2}{4r_0^2}}d\mu_0\\&&+\int_{r_0}^T\int_{\Sigma_{t-r_0^2}}|\bigtriangledown\phi
|| H |\frac{1}{4\pi
r_0^2}e^{-\frac{|F-X_0|^2}{4r_0^2}}d\mu_{t-r_0^2}dt\\&&+\int_{r_0}^T\int_{\Sigma_{t-r_0^2}}\frac{\phi}{8\pi
r_0^4}e^{-\frac{|F-X_0|^2}{4r_0^2}} |F-X_0||H| d\mu_{t-r_0^2}dt \\&&+ \int_{r_0}^T\int_{\Sigma_{t-r_0^2}}\frac{\phi}{4\pi
r_0^2}e^{-\frac{|F-X_0|^2}{4r_0^2}} |\tilde{R}-2| d\mu_{t-r_0^2}dt.
\end{eqnarray*} Using the $L^1$-estimate of the mean curvature vector and (\ref{e98}), note that $|\nabla\phi|\leq C$, we get that
\begin{eqnarray*}
\int_{\Sigma_{T-r_0^2}}\phi(F)\frac{1}{4\pi
r^2}e^{-\frac{|F-X_0|^2}{4r^2}}d\mu_{t-r^2} &\leq& 1+\varepsilon_0/2+\frac{C}{4\pi r_0^2}\frac{(\varepsilon_1+\varepsilon_2)^{\frac{1}{2}}}{1-e^{-\lambda/2}} \\&&
+\frac{C}{8\pi r_0^3}\frac{(\varepsilon_1+\varepsilon_2)^{\frac{1}{2}}}{1-e^{-\lambda/2}}+\frac{C\varepsilon_1}{4\pi r_0^2}.
\end{eqnarray*} If $(\varepsilon_1+\varepsilon_2)/r_0^6\ll \varepsilon_0$, then
\begin{eqnarray*}
\int_{\Sigma_{T-r_0^2}}\phi(F)\frac{1}{4\pi
r^2}e^{-\frac{|F-X_0|^2}{4r^2}}d\mu_{t-r^2} &\leq& 1+\varepsilon_0.
\end{eqnarray*} Applying White's theorem we obtain an uniform estimate of the second fundamental  from which
implies the global existence and convergence of the mean curvature flow. From (\ref{e96}) we see that
$\cos\alpha_\infty=1$, that is $\Sigma_\infty$ is a holomorphic curve. This proves the theorem.

\hfill Q. E. D

\end{document}